\documentclass[letterpaper, 10pt ,conference]{ieeeconf}

\IEEEoverridecommandlockouts  
\usepackage{graphicx}
\usepackage{dsfont} 
\usepackage{amsmath}
\usepackage{amssymb}
\usepackage{url}
\usepackage{algpseudocode}
\usepackage[utf8]{inputenc} 
\usepackage[format=hang,singlelinecheck=false,caption=false,font=footnotesize]{subfig}
\usepackage{color}

\graphicspath{{./Images/}}

\newtheorem{theorem}{Theorem}
\newtheorem{proposition}[theorem]{Proposition}

\newtheorem{assumption}[theorem]{Assumption}
\newtheorem{definition}[theorem]{Definition}

\let\originalleft\left
\let\originalright\right
\DeclareRobustCommand{\left}{\mathopen{}\mathclose\bgroup\originalleft}
\DeclareRobustCommand{\right}{\aftergroup\egroup\originalright}

\newcommand{\R}{\mathbb{R}}
\newcommand{\N}{\mathbb{N}}
\renewcommand{\epsilon}{\varepsilon}
\renewcommand{\phi}{\varphi}

\newcommand{\derivoper}{\mathrm{d}}

\DeclareMathOperator{\totvar}{V}
\DeclareMathOperator{\conv}{conv}

\newcommand{\norm}[1]{\left\lVert #1\right\rVert}
\newcommand{\normb}[1]{\bigl\lVert #1\bigr\rVert}

\newcommand{\set}[1]{\left\{ #1\right\}}

\newcommand{\setbb}[1]{\biggl\{ #1\biggr\}}

\let\p=\paren

\let\pb=\parenb

\let\pBB=\parenBB

\newcommand{\sqparen}[1]{\left[ #1\right]}
\newcommand{\sqparenb}[1]{\bigl[ #1\bigr]}

\newcommand{\dprod}[2]{\left\langle #1, #2\right\rangle}

\let\smat=\sqmatrix
\let\smats=\sqmatrixs

\newcommand{\diff}[1]{\derivoper #1}

\newcommand{\tp}{T}

\begin{document}

\title{%
  Optimal Control of Hybrid Systems Using a Feedback Relaxed Control Formulation%
}

\author{
  Tyler Westenbroek \and
  Humberto Gonzalez
  \thanks{%
    The authors are with the Department of Electrical \& Systems Engineering, Washington University in St.\ Louis, St.\ Louis, MO 63130.
    Email: {\scriptsize \texttt{\{westenbroek,hgonzale\}@wustl.edu}}.%
  }
}

\maketitle

\begin{abstract}
  We present a numerically tractable formulation for computing the optimal control of the class of hybrid dynamical systems whose trajectories are continuous.
  Our formulation, an extension of existing relaxed-control techniques for switched dynamical systems, incorporates the domain information of each discrete mode as part of the constraints in the optimization problem.
  Moreover, our numerical results are consistent with phenomena that are particular to hybrid systems, such as the creation of sliding trajectories between discrete modes.
\end{abstract}


\section{Introduction}
\label{sec:introduction}

Hybrid dynamical models are a powerful abstraction capable of describing the behavior of a wide range of systems such as legged locomotion robots~\cite{Ames2014}, unmanned aerial vehicles~\cite{Gillula2011}, or industrial process control plants~\cite{Gokbayrak2000,Li2015}, among many others.
Yet, even though hybrid models allow us to mathematically describe complex systems, efficient and flexible numerical methods to analyze and control those systems are hard to develop due to inherent phenomena such as Zeno executions~\cite{Johansson1999}.
In this paper we present a formulation for the class of hybrid systems whose trajectories are continuous (i.e., hybrid systems without reset maps) based on the theory of relaxed controls~\cite{Warga1972,Williamson1976,Warga1977}, and we show that our formulation allows us to numerically solve hybrid optimal control problems.

Before we can describe our results in details we must lay out the context in which our research is found.
Early results in optimality conditions for hybrid optimal control can be traced back to the work of Pontryagin et al.~\cite{Pontryagin1962}, who formulated a version of their Maximum Principle for systems whose vector field change discontinuously in predetermined stages.
In his seminal work, Sussmann~\cite{Sussmann1999,Sussmann1999a} formulated an extension of the Maximum Principle for a general class of hybrid systems incorporating autonomous and controlled discrete variables.
Yet, Sussmann's work is based upon a collection of sufficient conditions limiting the possible trajectories of the hybrid system in questions.
The Maximum Principle has been extended to other classes of hybrid systems, such as impulsive hybrid systems~\cite{Azhmyakov2003}.
Shaikh and Caines~\cite{Shaikh2003,Shaikh2007} also formulated an extension of the Maximum Principle for hybrid systems, and in their case they also found a numerical method capable of approximating the solution of those problems under some conditions.
Shaikh and Caines's optimality conditions impose restrictions to the optimal trajectories similar to Sussmann's optimality conditions.
In particular, the optimal trajectory cannot be Zeno (i.e., switching between discrete modes infinitely fast), nor non-orbitally stable (i.e., not continuous with respect to initial conditions), properties that cannot be easily checked in practice yet naturally occur in many applications~\cite{Lamperski2008}.
Hence, these restrictions have limited the deployment of numerical algorithms based on these optimality conditions to a wide range of applications.



The case of optimal control algorithms for switched hybrid systems deserves particular attention, since the properties of these systems have enabled the development of a number of numerical algorithms.
Xu and Antsaklis~\cite{Xu2002} formulated a general framework for the optimal control of switched systems dividing the problem in two parts, optimization of discrete sequences and optimization of discrete jump times.
Boccadoro et al.~\cite{Boccadoro2005} created a robust algorithm defining switching surfaces that are triggered by the continuous state variables, effectively creating a closed-loop solution.
Axelsson et al.~\cite{Axelsson2008} used needle variations to find an optimality condition that could be implemented using gradient descent algorithms.
Their approach was later extended by Gonzalez et al.~\cite{Gonzalez2010,Gonzalez2010a} to systems with constraints and continuous inputs.
Wardi and Egerstedt~\cite{Wardi2012a} used a measure theory approach to develop an switched optimal control algorithm.

Our work in this paper can be seen as an extension of the work by Morari et al.~\cite{Morari2003}, who transformed a discrete-time hybrid optimal control problem into a mixed integer optimization problem, Heemels and Brogliato~\cite{Heemels2003}, who used linear complementary problems to formulate discrete hybrid modes, and Vasudevan et al.~\cite{Vasudevan2013a,Vasudevan2013b}, who used a relaxed control formulation to create an adaptive refinement optimal control algorithm for switched systems.
In particular, we formulate the hybrid optimal control problem as a relaxed control switched optimal control problem (similar to~\cite{Vasudevan2013a}) with a set of extra constraints that use a complementary condition to disable certain discrete modes as a function of the continuous state values (similar to~\cite{Morari2003,Heemels2003}).
In practice, our formulation behaves as a switched hybrid system with a set-valued feedback law that takes the continuous state and returns a set of allowed discrete modes, which can be numerically approximated using standard techniques as a nonlinear programming problem~\cite{Polak1997}.

The paper is organized as follows.
Section~\ref{sec:math_bkgd} presents the mathematical background necessary to describe our results.
Section~\ref{sec:relaxation} presents our relaxed control formulation and the consequent hybrid optimal control problem.
Section~\ref{sec:num_imp} pre\-sents our numerical implementation of the hybrid optimal control.
Section~\ref{sec:examples} presents three hybrid optimal control problems solved using our algorithm.
Finally, Section~\ref{sec:conclusion} presents our conclusions and future research directions.

\section{Mathematical Background}
\label{sec:math_bkgd}

We begin by first formally defining the class of hybrid systems we will focus on.

\begin{definition}
  \label{def:hs}
  A controlled hybrid system with continuous trajectories is a tuple $\mathcal{H} = \pb{\mathcal{X}, \mathcal{I}, \mathcal{D}, \mathcal{U}, \mathcal{F}}$, where:
  \begin{itemize}
  \item $\mathcal{X} = \set{x \in \R^n \mid h(x) \leq 0}$ is the continuous-state space, with $h\colon \R^n \to \R^p$;
  \item $\mathcal{I} = \set{1,\dotsc,M}$ is the discrete-state space;
  \item $\mathcal{D} = \set{\mathcal{D}_i}_{i \in \mathcal{I}}$ is the set of domains, where for each $i \in {\cal I}$, $\mathcal{D}_i = \set{x \in {\cal X} \mid d_i(x) \leq 0}$ and $d_i\colon {\cal X} \to \R^{q_i}$;
  \item ${\cal U} \subset \R^m$ is the set control inputs; and,
  \item $\mathcal{F} = \set{f_i}_{i \in \mathcal{I}}$ is the set of vector fields, where $f_i\colon \mathcal{D}_i \times {\cal U} \to \R^n$ for each $i \in {\cal I}$.
\end{itemize}
\end{definition}

Throughout this paper we use $\norm{\cdot}_p$ to denote the finite-dimensional $p$-norm, i.e., given $x \in \R^n$, $\norm{x}_p = \sqrt[\leftroot{-3}\uproot{3}p]{\sum_{i=1}^n x_i^p}$.
We will abuse notation and denote the 2-norm simply by $\norm{\cdot}$.

\begin{definition}
  \label{def:l2}
  Let ${\cal S}$ be a subset of a Banach space and $\gamma\colon [0,T] \to {\cal S}$.
  We say that $\gamma \in L^2([0,T],{\cal S})$ iff:
  \begin{equation}
    \norm{\gamma}_{L^2} = \p{\int_0^T \norm{\gamma(t)}^2\, \diff{t}}^{\frac{1}{2}} < \infty.
  \end{equation}

  Also, let $\totvar(\gamma)$ be the \emph{total variation} of $\gamma$, defined by:
  \begin{multline}
    \totvar(\gamma) = \sup\setbb{\sum_{k=1}^N \normb{\gamma\p{t_{k+1}} - \gamma\p{t_k}}_1 \mid
      \set{t_k}_{k=1}^N \!\in {\cal P}_{[0,T]},\\[-5pt]
      N \in \N},
  \end{multline}
  where ${\cal P}_{[0,T]}$ is the set of all finite partitions of $[0,T]$.
  We say that $\gamma$ is of \emph{bounded variation}, denoted $\gamma \in BV([0,T],{\cal S})$, iff $\totvar(\gamma) < \infty$.
\end{definition}

We impose the following standard assumptions to ensure the existence and uniqueness of the continuous trajectories in each discrete mode, together with the well-posedness of the constraints in our optimization problem.
Details about the practical consequences of these assumptions can be found in Chapter~5 of~\cite{Polak1997}.

\begin{assumption}
  \label{assump:lipschitz}
  The functions $\set{f_i}_{i \in {\cal I}}$, $h$, and $\set{d_i}_{i \in {\cal I}}$ are continuously differentiable.
  Moreover, all these functions are globally Lipschitz continuous, and so are all their partial derivatives.
\end{assumption}

\begin{assumption}
  \label{assump:qc}
  The set ${\cal U}$ is compact.
  Also, the sets ${\cal U}$, ${\cal X}$, and $\set{{\cal D}_i}_{i \in {\cal I}}$ are constraint qualified.
\end{assumption}
More information regarding constraint qualified sets can be found in Chapter~12.2 of~\cite{Nocedal2006}.

Now we state our definition of hybrid execution and hybrid trajectory.
\begin{definition}
  \label{def:exec}
  Let ${\cal H}$ be a hybrid system as in Definition~\ref{def:hs}, $u \in L^2([0,T],{\cal U})$, and $x\colon [0,T] \to {\cal X}$.
  We say that the pair $(x,u)$ is a hybrid execution of ${\cal H}$ with length $T$ iff for almost every $t \in [0,T]$ there exists $i \in {\cal I}$ such that:
  \begin{align}
    \label{eq:exec1}
    &\dot{x}(t) = f_i\pb{x(t), u(t)}; \quad \text{and},\\
    \label{eq:exec2}
    &d_i\pb{x(t)} \leq 0.
  \end{align}
  Under these conditions we also say that $x$ is a trajectory of ${\cal H}$.
\end{definition}
Note that Definition~\ref{def:exec} implies that every trajectory of a hybrid system is absolutely continuous (as defined, for example, in~\cite{Folland1999}), which is consistent with the class of hybrid system in Definition~\ref{def:hs}.
Also note that we make no assumptions regarding domains being disjoint, executions being unique, or Zeno phenomena not appearing.
On the contrary, as we show in Section~\ref{sec:examples} we can naturally formulate problems where domains overlap (therefore executions are not unique), as well as problems where the optimal solutions are either Zeno or sliding modes (as understood in the Filippov sense~\cite{Filippov1988}).

As explained in Section~\ref{sec:introduction}, our main result extends the relaxed control formulation for switched systems in~\cite{Vasudevan2013a} to hybrid systems with continuous trajectories.
Hence, using the notation in~\cite{Vasudevan2013a}, we introduce the $M$-simplex and the set of corners of the $M$-simplex.
\begin{definition}
  \label{def:simplex}
  Let $M \in \N$.
  We define the \emph{$M$-simplex}, denoted $\Sigma_r^M$, and the set of \emph{corners of the $M$-simplex}, denoted $\Sigma_p^M$, by:
  \begin{align}
    \Sigma_r^M &= \setbb{\omega \in \R^M \mid \sum_{i=1}^M \omega_i = 1,\ \omega_i \in [0,1]\ \forall i};\ \text{and},\\
    \Sigma_p^M &= \setbb{\omega \in \R^M \mid \sum_{i=1}^M \omega_i = 1,\ \omega_i \in \set{0,1}\ \forall i}.
  \end{align}
\end{definition}

Note that $\Sigma_p^M$ contains $M$ elements, thus if $\omega \in \Sigma_p^M$ then exactly one of its entries $\omega_i = 1$ while all the rest are zero.
Thanks to this property we can use vectors in $\Sigma_p^M$ at each time $t$ to indicate the active discrete mode of a hybrid system.
Also note that $\Sigma_r^M$ is in fact the space of measures (or, depending of the context, probability distributions) defined over the set $\Sigma_p^M$.
Hence, if we choose a vector in $\Sigma_r^M$ to indicate the current discrete mode we are in fact choosing a relaxed control representation (as described in~\cite{Warga1972}) of a vector in $\Sigma_p^M$.
This property inspires our notation, where the subscript $p$ stands for \emph{pure} indicator vectors, while the subscript $r$ stands for \emph{relaxed} indicator vectors.

We are now ready to present our main result involving the use of a relaxed control formulation to approximate hybrid optimal control problems.

\section{Relaxed Control Formulation}
\label{sec:relaxation}

Our relaxed control formulation is based on two transformations, one for~\eqref{eq:exec1} and one for~\eqref{eq:exec2}, aimed at enabling a simple numerical formulation to compute hybrid executions.
First, note that if there exists $\omega \in L^2([0,T], \Sigma_p^M)$ such that:
\begin{equation}
  \label{eq:rexec1}
  \dot{x}(t) = \sum_{i=1}^M \omega_i(t)\, f_i\pb{x(t),u(t)},
\end{equation}
then, since exactly one entry $\omega_i(t) = 1$ for almost every $t \in [0,T]$, we recover the condition in~\eqref{eq:exec1}.
In other words, a function $\omega \in L^2([0,T],\Sigma_p^M)$ behaves as an indicator of the active discrete mode of the hybrid system at time $t$.
As shown in Theorem~3.3 of~\cite{Vasudevan2013a}, the differential equation in~\eqref{eq:rexec1} is well defined with unique solutions in the interval $[0,T]$.

Second, note that given $\omega \in L^2([0,T],\Sigma_p^M)$, the condition in~\eqref{eq:exec2} implies that if $d_i\pb{x(t)} > 0$ then $\omega_i(t) = 0$, i.e., if $x(t)$ does not belong to ${\cal D}_i$ then indicator vector cannot set mode $i$ as active.
The contrapositive is therefore also true, i.e., if $\omega_i(t) = 1$ then $d_i\pb{x(t)} \leq 0$.
It is worth noting that the condition in~\eqref{eq:exec2} does not impose restrictions to $\omega$ when $d_i\pb{x(t)} \leq 0$, since in that case any other mode $j \in {\cal I}$ might be active at time $t$ provided $d_j\pb{x(t)} \leq 0$.
We summarize the relation between $\omega_i(t)$ and $d_i\pb{x(t)}$ as follows:
\begin{equation}
  \label{eq:rexec2}
  \omega_i(t)\, d_i\pb{x(t)} \leq 0, \quad \forall i \in {\cal I}.
\end{equation}
In this way, when each of the inequalities in~\eqref{eq:rexec2} is satisfied we get a complementarity condition where both $\omega_i(t) = 1$ and $x(t) \notin {\cal D}_i$ cannot occur at the same time $t$.

The enforcement of the conditions in equations~\eqref{eq:rexec1} and~\eqref{eq:rexec2} can be theoretically interpreted as a switched hybrid system with a feedback rule disabling certain discrete modes as a function of $x(t)$.
Indeed, let:
\begin{equation}
  {\cal S}_p\pb{x(t)} = \set{\omega \in \Sigma_p^M \mid \omega_i = 0\ \text{whenever}\ d_i\pb{x(t)} > 0},
\end{equation}
then the condition in~\eqref{eq:rexec2} is equivalent to enforcing $\omega(t) \in {\cal S}_p\pb{x(t)}$.
This interpretation opens the doors to use algorithms computing the optimal control of switched systems in a bigger class of hybrid systems, as the ones in Definition~\ref{def:hs}.
In particular, we will use the approach defined in~\cite{Vasudevan2013a} to compute the optimal control of hybrid systems.

The proof of the following proposition follows directly from the argument described above.
\begin{proposition}
  \label{prop:exec}
  Let ${\cal H}$ be a hybrid system as in Definition~\ref{def:hs}, $u \in L^2([0,T],{\cal U})$, and $x\colon [0,T] \to {\cal X}$.

  Then $(x,u)$ is a hybrid execution of ${\cal H}$ with length $T$ iff there exists $\omega \in L^2([0,T],\Sigma_p^M)$ such that~\eqref{eq:rexec1} and~\eqref{eq:rexec2} are satisfied.
\end{proposition}

Using on Proposition~\ref{prop:exec} we get the following optimal control problem.
\begin{definition}
  \label{def:hocp}
  Let ${\cal H}$ be a hybrid system as in Definition~\ref{def:hs}.
  A \emph{hybrid optimal control problem} is defined by:
  \begin{equation}
    \label{eq:hocp}
    \begin{aligned}
      \min_{x(t), u(t), \omega(t)} &\int_0^T \sum_{i=1}^M \omega_i(t)\, L_i\pb{x(t),u(t)}\, \diff{t} + \phi\pb{x(T)},\\
      \text{s.t.}\
      &x(0) = \xi,\\
      &\dot{x}(t) = \sum_{i=1}^M \omega_i(t)\, f_i\pb{x(t),u(t)},\ \forall t \in [0,T],\\
      &\omega_i(t)\, d_i\pb{x(t)} \leq 0,\ \forall i \in {\cal I},\ \forall t \in [0,T],\\
      &h\pb{x(t)} \leq 0,\ \forall t \in [0,T],\\
      &u(t) \in {\cal U},\ \forall t \in [0,T],\\
      &\omega(t) \in \Sigma_p^M,\ \forall t \in [0,T].
    \end{aligned}
  \end{equation}
  where $\set{L_i}_{i=1}^M$ and $\phi$ are continuously differentiable real-valued functions, and $\xi \in \R^n$ is the initial condition of the continuous-state variables.
\end{definition}
Note that since $\Sigma_p^M$ is a discrete set with exactly $M$ elements, the optimal control problem in Definition~\ref{def:hocp} is in practice a mixed-integer programming problem, which cannot be efficiently solved using gradient-based numerical algorithms.

\subsection{Optimal Control using Relaxed Controls}

In view of the result in Proposition~\ref{prop:exec} and the relation between $\Sigma_p^M$ and $\Sigma_r^M$, we can extend Definition~\ref{def:exec} to consider \emph{relaxed executions}.
\begin{definition}
  \label{def:rexec}
  Let ${\cal H}$ be a hybrid system as in Definition~\ref{def:hs}, $u \in L^2([0,T],{\cal U})$, and $x\colon [0,T] \to {\cal X}$.

  Then $(x,u)$ is a relaxed hybrid execution of ${\cal H}$ with length $T$ iff there exists $\omega \in L^2([0,T],\Sigma_r^M)$ such that~\eqref{eq:rexec1} and~\eqref{eq:rexec2} are satisfied.
\end{definition}

We can now avoid formulating the mixed-integer program in~\eqref{eq:hocp} by modifying its last constraint.
\begin{definition}
  \label{def:rhocp}
  Let ${\cal H}$ be a hybrid system as in Definition~\ref{def:hs}.
  A \emph{relaxed hybrid optimal control problem} is defined as a relaxation of the hybrid optimal control problem in Definition~\ref{def:hocp}, where we replace the last constraint $\omega(t) \in \Sigma_p^M$ with $\omega(t) \in \Sigma_r^M$ in~\eqref{eq:hocp}.
\end{definition}

Recall that $\Sigma_p^M \subset \Sigma_r^M$, hence the feasible set of the relaxed problem is larger than that of the pure, and should therefore result in a lower optimal value in theory.
Yet as shown in~\cite{Vasudevan2013a}, and explained below, these two optimization problems produce exactly the same optimal value.

Consider the following operators:
\begin{definition}
  \label{def:haar}
  Let $b_{0,0}\colon [0,T] \to \R$ be the \emph{Haar wavelet}, defined by:
  \begin{equation}
    b_{0,0}(t) =
    \begin{cases}
      1 & \text{if}\ t \leq \frac{T}{2},\\
      0 & \text{if}\ t > \frac{T}{2}.
    \end{cases}
  \end{equation}
  Also, for each $k \in \N$ and $j \in \set{0,\dotsc,2^k - 1}$ let $b_{k,j}(t) = b_{0,0}\pb{T\,\p{2^k\, t - j}}$.

  We say that ${\cal F}_N\colon L^2([0,T],\Sigma_r^M) \to L^2([0,T],\Sigma_r^M)$ is the \emph{Haar wavelet operator}, defined by:
  \begin{equation}
    \sqparenb{{\cal F}_N(\omega)}_i(t) = \frac{1}{T}\, \dprod{\omega_i}{\mathds{1}} + \sum_{k=0}^N \sum_{j=0}^{2^k-1} \dprod{\omega_i}{b_{k,j}}\, \frac{b_{k,j}(t)}{\norm{b_{k,j}}_{L^2}^2},
  \end{equation}
  for each $i \in {\cal I}$.
\end{definition}
In other words, given $\omega \in L^2([0,T],\Sigma_r^M)$ the Haar wavelet operator returns a new function ${\cal F}_N(\omega) \in L^2([0,T],\Sigma_r^M)$ that is piecewise constant over a uniform partition of $2^N$ samples.

\begin{definition}
  \label{def:pwm}
  We define the \emph{pulse-width modulation operator}, denoted ${\cal P}_N\colon L^2([0,T],\Sigma_r^M) \to L^2([0,T],\Sigma_p^M)$, by:
  \begin{equation}
    \sqparenb{{\cal P}_N(\omega)}_i(t) =
    \begin{cases}
      1 & \text{if}\ t \in \left[ S_{k,i-1}, S_{k,i} \right),\\
      0 & \text{otherwise},
    \end{cases}
  \end{equation}
  where for each $k \in \set{0,\dotsc,2^N\!\!-1}$ and $i \in \set{1,\dotsc,M}$:
  \begin{equation}
    S_{k,i} = \frac{T}{2^N}\, \pBB{k + \sum_{j=1}^{i} \omega_j\p{\frac{k}{2^N}}}.
  \end{equation}
\end{definition}
Hence, given $\omega \in L^2([0,T],\Sigma_r^M)$ the pulse-width modulation operator returns a new function ${\cal P}_N(\omega) \in L^2([0,T],\Sigma_p^M)$ where each coordinate switches between $0$ and $1$ with pulses whose width is proportional to the amplitude of $\omega$ evaluated at the samples $\frac{k}{2^N}$ for $k \in \set{0,\dotsc,2^N\!\!-1}$.

Using these two operators we get the following result.
The proof is omitted since it is an extension of Theorem~5.10 in~\cite{Vasudevan2013a} for the case when $\omega(t) \in \conv\pb{{\cal S}_p\pb{x(t)}}$, the convex hull of ${\cal S}_p\pb{x(t)}$, for each $t \in [0,T]$.
\begin{theorem}
  \label{thm:approx}
  Let ${\cal H}$ be a hybrid system as in Definition~\ref{def:hs}, and let $\omega \in L^2([0,T],\Sigma_r^M)$.

  Furthermore, let $(x,u)$ be the hybrid execution associated to $\omega$, and for each $N \in \N$, let $(x_N,u)$ be the hybrid execution associated to ${\cal P}_N\pb{{\cal F}_N(\omega)} \in L^2([0,T],\Sigma_p^M)$.

  If $\omega \in BV([0,T],\Sigma_r^M)$, then:
  \begin{equation}
    \normb{x(t) - x_N(t)} \leq C\, 2^{-\frac{N}{2}}\, \pb{\totvar(\omega) + 1},
  \end{equation}
  where $C > 0$ does not depend on $\omega$.
\end{theorem}

Theorem~\ref{thm:approx} implies that, provided the optimal indexing function $\omega$ of the relaxed optimal control problem in Definition~\ref{def:rhocp} is of bounded variation, then we can create a sequence of points in the feasible set of the optimization problem in Definition~\ref{def:hocp} that approximates the optimal relaxed execution.
This fact implies that, in most situations, both optimization problems produce exactly the same value, and therefore by solving the relaxed hybrid optimal control problem we are in fact solving the original hybrid optimal control problem as well.
It is worth noting that both projection operators, ${\cal F}_N$ and ${\cal P}_N$, can be numerically implemented using efficient algorithms (an implementation is provided in~\cite{optwrapper}).

A discussion regarding whether the optimal indexing function $\omega$ of the relaxed hybrid optimal control problem is of bounded variation is beyond the scope of this paper.
For the time being it is worth mentioning that all functions with weak derivatives (as defined in~\cite{Ziemer1989}) belong to the set of bounded variation, and that all functions belonging to a finite-dimensional subspace of $L^2$, such as the functions obtained using numerical methods, also belong to the set of bounded variation.

At this point we have shown that the optimal relaxed solutions resulting from the relaxed hybrid optimal control problem approximate solutions to the original hybrid optimal control problem with arbitrary accuracy.
In the next section we detail the numerical implementation of the results above.

\section{Numerical Implementation}
\label{sec:num_imp}

Our numerical implementation of the optimization problem in Definition~\ref{def:rhocp} uses the Forward Euler discretization to approximate the differential equation constraints and the integral term in the cost function.
Even though other algorithms such as Runge-Kutta or Pseudospectral methods are more accurate than the Forward Euler discretization~\cite{Gong2006}, we chose to use the latter due to its simple implementation and good convergence properties when used in optimal control problems (as shown in Chapter~4 of~\cite{Polak1997}).




Applying the Forward Euler discretization to our formulation, we get the following nonlinear programming optimization problem, which approximates the problem in Definition~\ref{def:rhocp}:
\begin{equation}
  \label{eq:disc_rhocp}
  \begin{aligned}
    \min_{\set{\bar{x}_k}, \set{\bar{u}_k}, \set{\bar{\omega}_k}}
    &\sum_{k=0}^{N-1} \sum_{i=1}^M \frac{T}{N}\, \bar{\omega}_{k,i}\, L_i\p{x_k,u_k} + \phi\pb{x_N},\\
    \text{s.t.}\
    &\bar{x}_0 = \xi,\\
    &\bar{x}_{k+1} \!\!=\! \bar{x}_k \!+\! \frac{T}{N} \!\sum_{i=1}^M \bar{\omega}_{k,i}\, f_i\p{\bar{x}_k,\bar{u}_k}, \forall k \!\leq \!N\!-\!1,\\
    &\bar{\omega}_i\, d_i\p{\bar{x}_k} \leq 0,\ \forall i \in {\cal I},\ \forall k \leq N,\\
    &h\p{\bar{x}_k} \leq 0,\ \forall k \leq N,\\
    &\bar{u}_k \in {\cal U},\ \forall k \leq N-1,\\
    &\bar{\omega}_k \in \Sigma_r^M,\ \forall k \leq N.
  \end{aligned}
\end{equation}

The nonlinear programming problem in~\eqref{eq:disc_rhocp} results in a relaxed hybrid execution $(x,u)$ associated to the relaxed indexing mode function $\omega$.
Then, given $N \in \N$ we can find a hybrid execution $(x_N,u)$ associated to the indexing function ${\cal P}_N\pb{{\cal F}_N(\omega)}$, where ${\cal F}_N$ and ${\cal P}_N$ are as in Definitions~\ref{def:haar} and~\ref{def:pwm}, respectively.

In the next section we present two examples where we apply this numerical implementation to solve three hybrid optimal control problems.

\section{Examples}
\label{sec:examples}

\subsection{Sliding-Mode Trajectory}
\label{sec:ex_sm}

In this example we show how our formulation produces optimal sliding-mode trajectories.
Consider the hybrid system with continuous-state space ${\cal X} = \R^2$, no control inputs (i.e., ${\cal U} = \emptyset$), and with two overlapping domains as defined in Table~\ref{tab:sm_doms}.
Note that this hybrid system does not have unique trajectories, since any of the modes can be chosen in the intersection between their domains.

\begin{table}[tp]
  \centering
  \caption{%
    Discrete modes of the hybrid system in Example~\ref{sec:ex_sm}.%
  }%
  \label{tab:sm_doms}%
  \begin{tabular}{c|c|c}
    Mode  & Domain  & Vector Field \\
    \hline
    1 & $x_2 \geq \sin(x_1) - \frac{1}{2}$ & $\dot{x} = \smat{1 & -1}^\tp$ \rule{0pt}{2.6ex}\\[1pt]
    2 & $x_2 \leq \sin(x_1) + \frac{1}{2}$ & $\dot{x} = \smat{1 & 1}^\tp$
  \end{tabular}
\end{table}

Using the notation in~\eqref{eq:disc_rhocp}, the cost function is defined by $L_1(x,u) = L_2(x,u) = \phi(x) = \norm{x}^2$, the initial condition is $x_0 = \smats{0 \\ -1}$, the optimization horizon is $T=4.5$, and the discretization has $N=100$ samples.
The results were obtained using the \emph{SNOPT} solver~\cite{Gill2002}, and the calculations took $6.28$~seconds in an \emph{Intel Xeon E5-2680} processor running at $2.7$~GHz.

\begin{figure*}[tp]
  \centering
  \subfloat[%
  Optimal relaxed trajectory.
  Top region (green) is the domain of Mode~1, bottom region (red) is the domain of Mode~2, and middle region (white) is the overlap of both Modes.%
  ]{%
    \label{fig:ex1a}%
    \includegraphics[width=.32\linewidth]{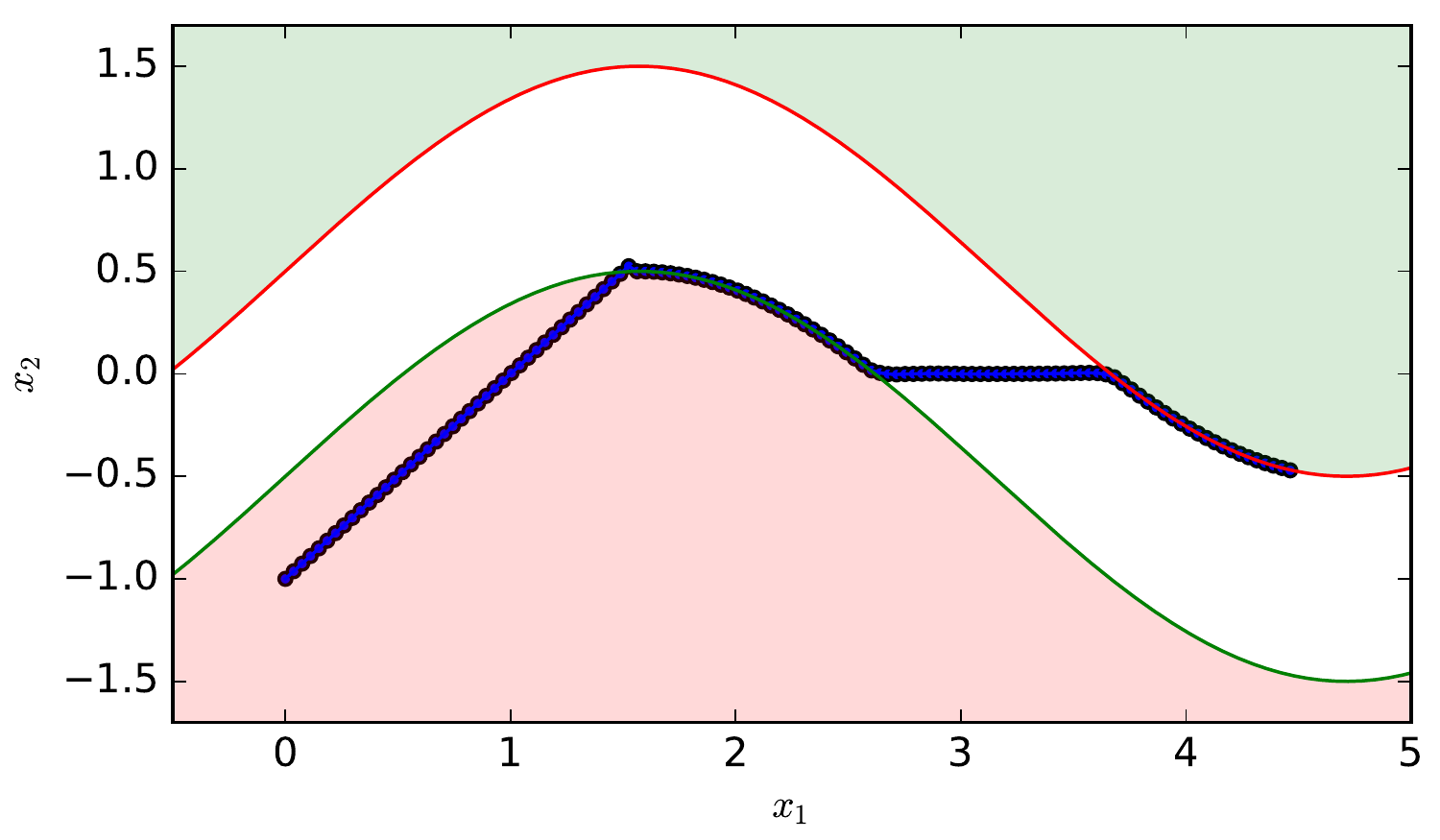}%
  }%
  \hfill%
  \subfloat[%
  Hybrid trajectory using the projection operator ${\cal P}_N\pb{{\cal F}_N(\omega)}$, with $N=5$, on the optimal solution in Fig.~\ref{fig:ex1a}.%
  ]{%
    \label{fig:ex1c}%
    \includegraphics[width=.32\linewidth]{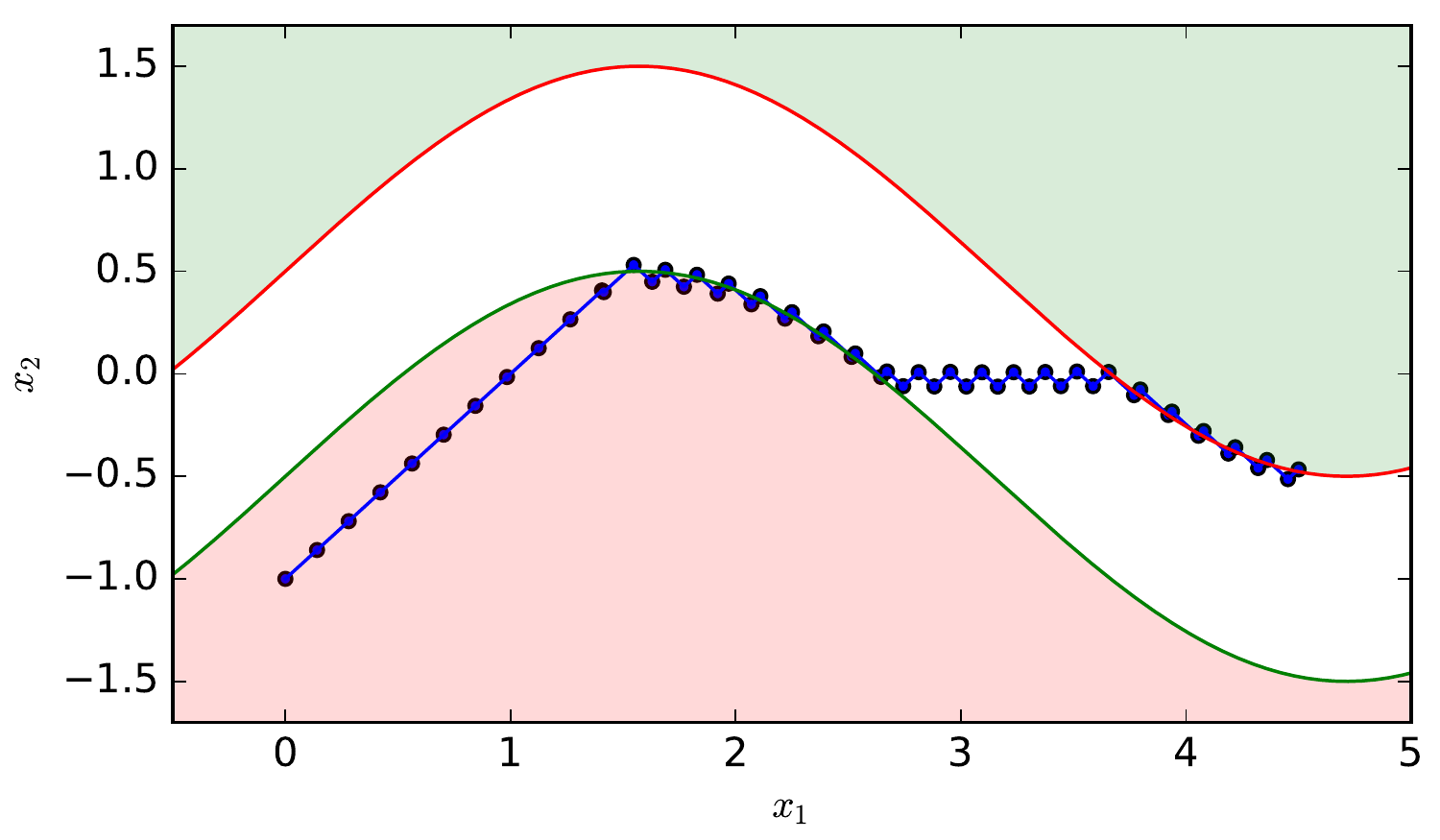}%
  }%
  \hfill%
  \subfloat[%
  Indexing function $\omega_1$ associated to Mode~1.
  \emph{Top}: optimal relaxed indexing function.
  \emph{Bottom}: projected indexing function with $N=5$.
  ]{%
    \label{fig:ex1b}%
    \includegraphics[width=.32\linewidth]{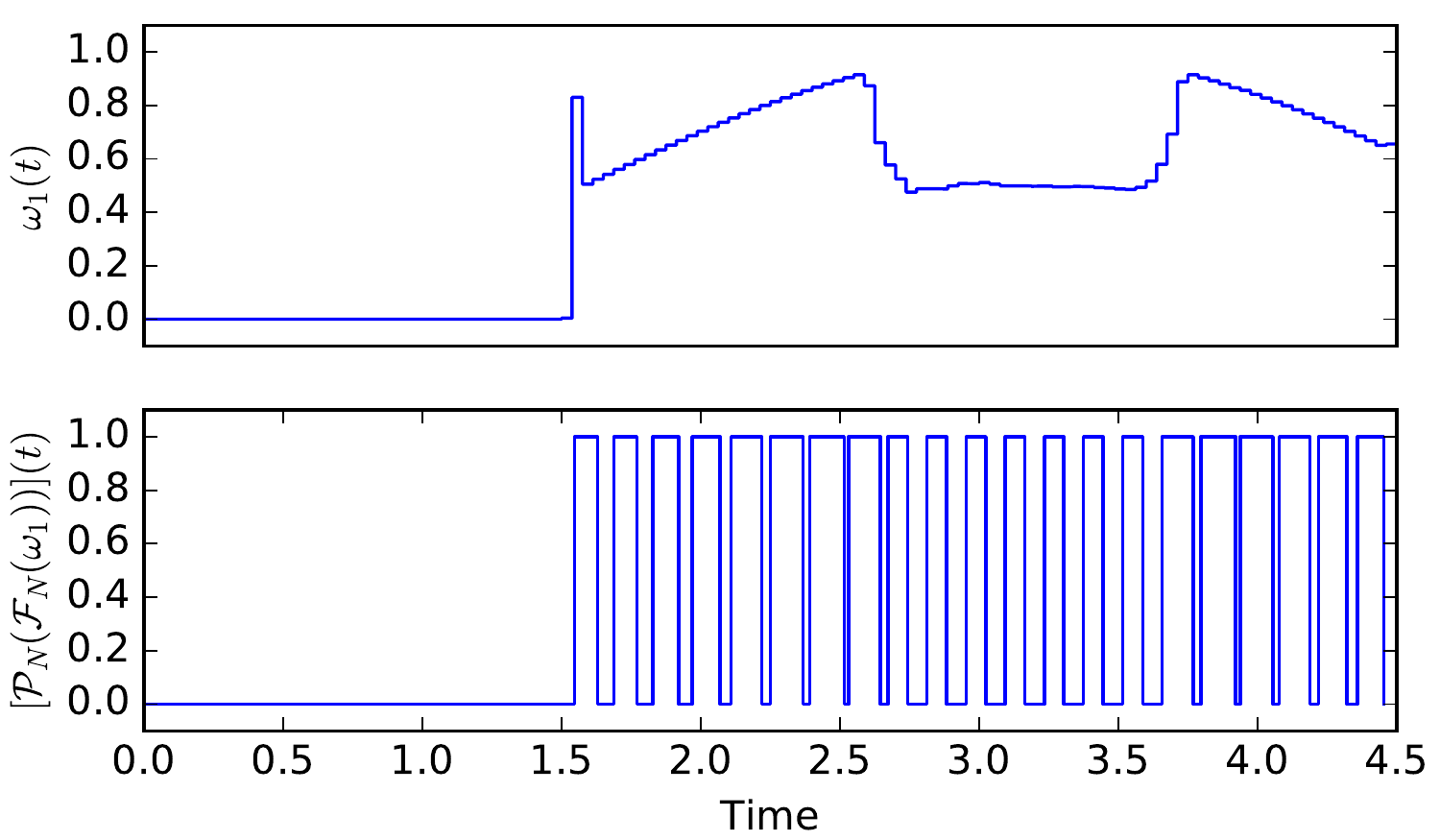}%
  }%
  \caption{Results of Example~\ref{sec:ex_sm}.}
  \label{fig:ex1}
\end{figure*}

As shown in Figure~\ref{fig:ex1}, the optimal solution of this problem is a sliding-mode trajectory, where modes~1 and~2 chatter in the intersection between their domains.
Thus, as shown in Figure~\ref{fig:ex1c}, the relaxed indexing function $\omega(t)$ contains non-binary values, i.e., in the relaxed trajectory both discrete modes are applied simultaneously.
Note that the execution generates two different types of sliding modes.
In the first type, the trajectory follows the boundary between the two modes, and the modal weights converge to a sinusoidal pattern in order to keep the trajectory along this surface. In the second type of sliding mode, $\omega_1(t) \approx \omega_2(t) \approx \frac{1}{2}$ and the trajectory moves horizontally through the interior of the intersection of the two modes.
This ability to transition between different sliding modes, under the influence of the objective function, allows our formulation to pick a relaxed optimal trajectory through regions of the state-space where Zeno phenomena readily occur.

Figures~\ref{fig:ex1b} and~\ref{fig:ex1c} show the effect of using the operators ${\cal F}_N$ and ${\cal P}_N$ to generate new hybrid indexing functions that, although sub-optimal, follow closely the behavior of the optimal relaxed solution, even for small values of $N$.

\subsection{Partitioned Linear System}
\label{sec:ex_lin}

In this example we again consider ${\cal X} = \R^2$, but now we partition the plane along the lines $x_1 = x_2$ and $x_1 = - x_2$, creating four distinct modes with no overlap.
The domains and vector fields are summarized in Table~\ref{tab:lin_doms}, where non-adjacent modes are endowed with the same vector field.
The vector field matrices are $A_1 = \smats{-0.1 & 0.1 \\ -0.4 & -0.1}$, and $A_2 = \smats{0.1 & 0.5 \\ -0.3 & -0.1}$.
The control input set is defined by ${\cal U} = \set{u \in \R \mid u \in [-1,0.1]}$, and the cost function is defined by $L(x,u) = \norm{x}^2 + \norm{u}^2$ for all modes, and $\phi(x) = \norm{x}^2$.

\begin{figure}[tp]
  \centering
  \includegraphics[width=.66\linewidth]{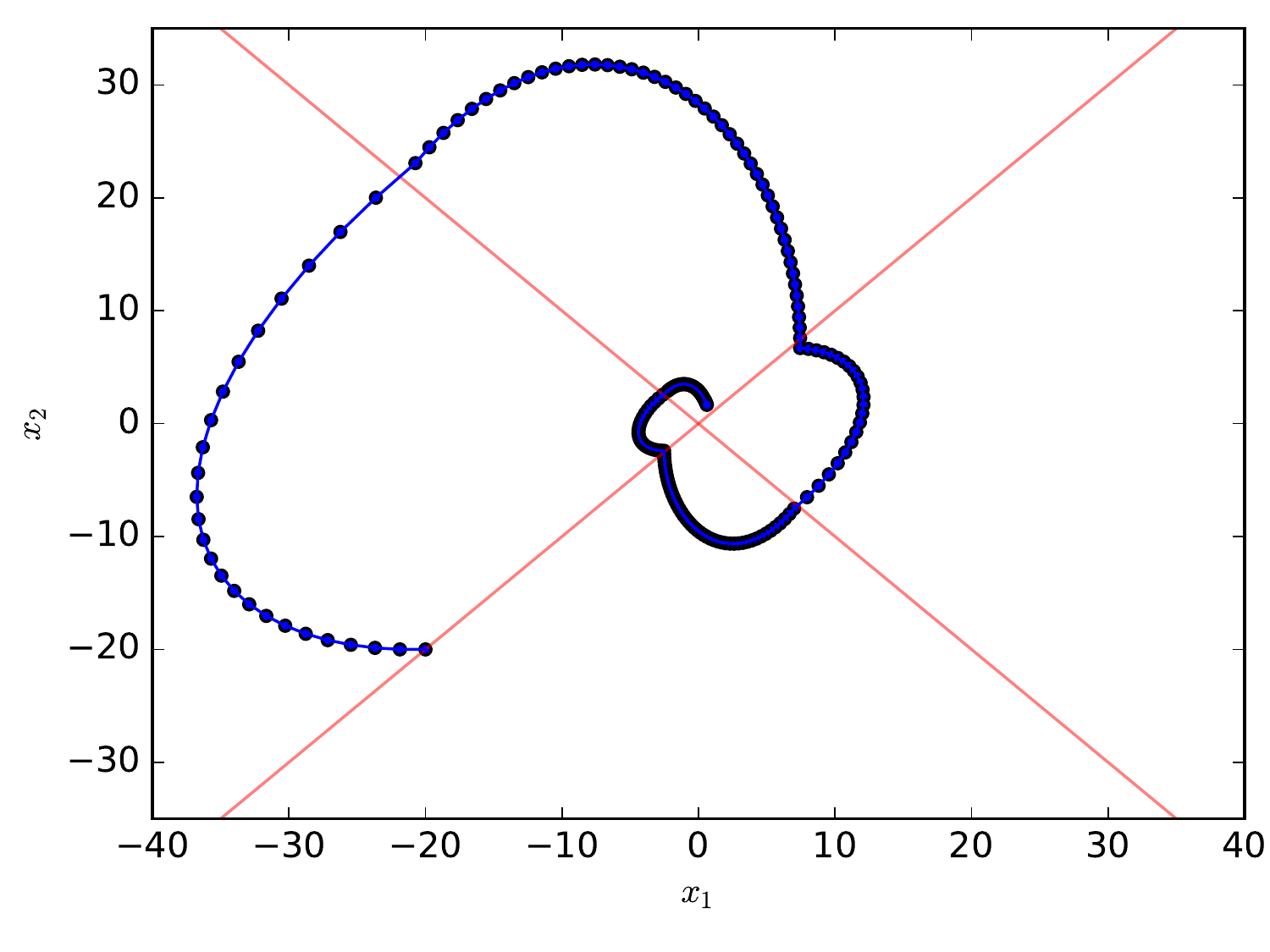}%
  \caption{%
    Optimal trajectory of Example~\ref{sec:ex_lin}.
    The boundaries between modes, as defined in Table~\ref{tab:lin_doms}, are shown in red.
  }%
  \label{fig:ex2}
\end{figure}

Figure~\ref{fig:ex2} shows the optimal result with initial condition $x_0 = \smats{0 \\ -1}$, horizon $T=40$, and $N=100$ samples.
Note that hybrid systems without reset maps and non-overlapping domains are commonly used in control applications due to their easy design and intuitive formulation~\cite{Branicky1998a}.
Our result shows that the problem in~\eqref{eq:disc_rhocp} can be used to efficiently solve this large class of hybrid systems.

\begin{table}[tp]
  \centering
  \caption{%
    Discrete modes of the hybrid system in Example~\ref{sec:ex_lin}.
  }%
  \label{tab:lin_doms}%
  \begin{tabular}{c|c|c}
    Mode  & Domain  & Vector Field \\
    \hline
    1 & $-x_2 \leq x_1 \leq x_2$ & $\dot{x} = A_1\, x(t) + B\, u(t)$ \rule{0pt}{2ex}\\
    2 & $x_1 \leq x_2 \leq -x_1$ & $\dot{x} = A_2\, x(t) + B\, u(t)$\\[1pt]
    3 & $x_2 \leq x_1 \leq -x_2$ & $\dot{x} = A_1\, x(t) + B\, u(t)$\\[1pt]
    4 & $-x_1 \leq x_2 \leq x_1$ & $\dot{x} = A_2\, x(t) + B\, u(t)$
  \end{tabular}
\end{table}

\subsection{Quadrotor Helicopter Obstacle Avoidance}
\label{sec:ex_qr}

We now consider a simplified quadrotor helicopter model with six continuous states, similar to the one used in Section~5.3 of~\cite{Vasudevan2013b}.
The first state is the horizontal displacement of the helicopter, the second state is its vertical displacement, and the third state is its the roll angle.
The remaining three states are the respective derivatives of the first three states.
We model the helicopter as a hybrid system with three discrete modes: one whose task is to move the helicopter vertically (denoted $U$), the second to rotate to the left (denoted $L$), and the third to rotate to the right (denoted $R$).
Each of these discrete modes has a single input.

Our control objective is to make the helicopter lift off the ground, fly over an obstacle, and land safely at the origin.
Even though we could incorporate the obstacle as part of the domain constraints in~\eqref{eq:disc_rhocp} (i.e., using a nonlinear constraint function $h$), that approach typically results in longer computation times.
Instead, we use a penalty function to increase the cost of flying through a region around the obstacle.
We also add a fourth safety mode (denoted $S$), in which the helicopter is permitted to accelerate upwards very quickly, at the expense of higher battery usage, in order to pass over the obstacle.
The input is not penalized in the mode, in order to incentivize the use of this extra mode when the helicopter is at risk of collision.

The discrete modes in this example are detailed in Table~\ref{tab:qr_doms}, with the following parameters: gravity constant $g=9.8$, helicopter mass $\mu=1.3$, horizontal span $\ell=0.305$, moment of inertia $I=0.0605$, and motor constant $\lambda = 1530$.

\begin{table}[tp]
  \centering
  \caption{%
    Discrete modes of the hybrid system in Example~\ref{sec:ex_qr}.%
  }%
  \label{tab:qr_doms}%
  \begin{tabular}{c|c|c}
    Mode  & Domain  & Vector Field \\
    \hline
    $U$ & $ x_3 \in [-0.1,0.1] $ & $\ddot{x}(t) = \smat{\frac{\lambda}{\mu} \sin(x_3)\, u\\ \frac{\lambda}{\mu} \cos(x_3)\, u - g\\ 0}$  \rule{0pt}{5.8ex}\\[14pt]
    $L$ & $ x_3 \geq -0.1$ & $\ddot{x}(t) = \smat{g \sin(x_3)\\ g \cos(x_3) - g\\ - \frac{\ell}{I}\, u}$\\[12pt]
    $R$ & $ x_3 \leq 0.1$  & $\ddot{x}(t) = \smat{g \sin(x_3)\\ g \cos(x_3) - g\\ \frac{\ell}{I}\, u}$\\[12pt]
    $S$ & $2\,(x_1-4)^2 + x_2^2 \leq 4$ & $\ddot{x}(t) = \smat{\frac{\lambda}{\mu} \sin(x_3)\, u\\ \frac{3\, \lambda}{2\, \mu} \cos(x_3)\, u - g\\ 0}$\\
  \end{tabular}
\end{table}

Our cost function is defined by $L_U(x,u) = L_R(x,u) = L_L(x,u) = 5\, u^2 + \rho(x)$, $L_S(x,u) = \rho(x)$, and $\phi(x) = 10^3\, \p{x_1^2 + x_2^2 + x_4^2} + 10^2\, x_3^2$, where $\rho$ is the \emph{bump} function defined by $\rho(x) =\exp\p{\frac{-1}{1 - 2\, \p{x_1-4}^2 - x_2^2}}$, if $2(x_1-4)^2 + x_2^2 \leq 1$, and $\rho(x)=0$ otherwise.
Furthermore, we constraint the input to $u(t) \in [0,1]$. 
We also add the safety constraint $x_2(t) \geq 0$ so that the helicopter stays above ground.
The time horizon is $T=6$, with $N=100$ samples, and the initial condition is $x_1(0) = 6$ with all other states set equal to zero.
The simulation, also solved using \emph{SNOPT} in the same computer as the previous examples, took~201 seconds to complete.

\begin{figure*}[tp]
  \centering
  \subfloat[%
  Optimal relaxed trajectory.
  The red ellipse shows the outline of the obstacle, the green ellipse shows the domain of Mode~$S$.
  Red dots mark the use of Mode~$S$.
  A sketch of the helicopter is shown in black for a handful of time samples.%
  ]{%
    \label{fig:ex3a}%
    \includegraphics[width=.32\linewidth]{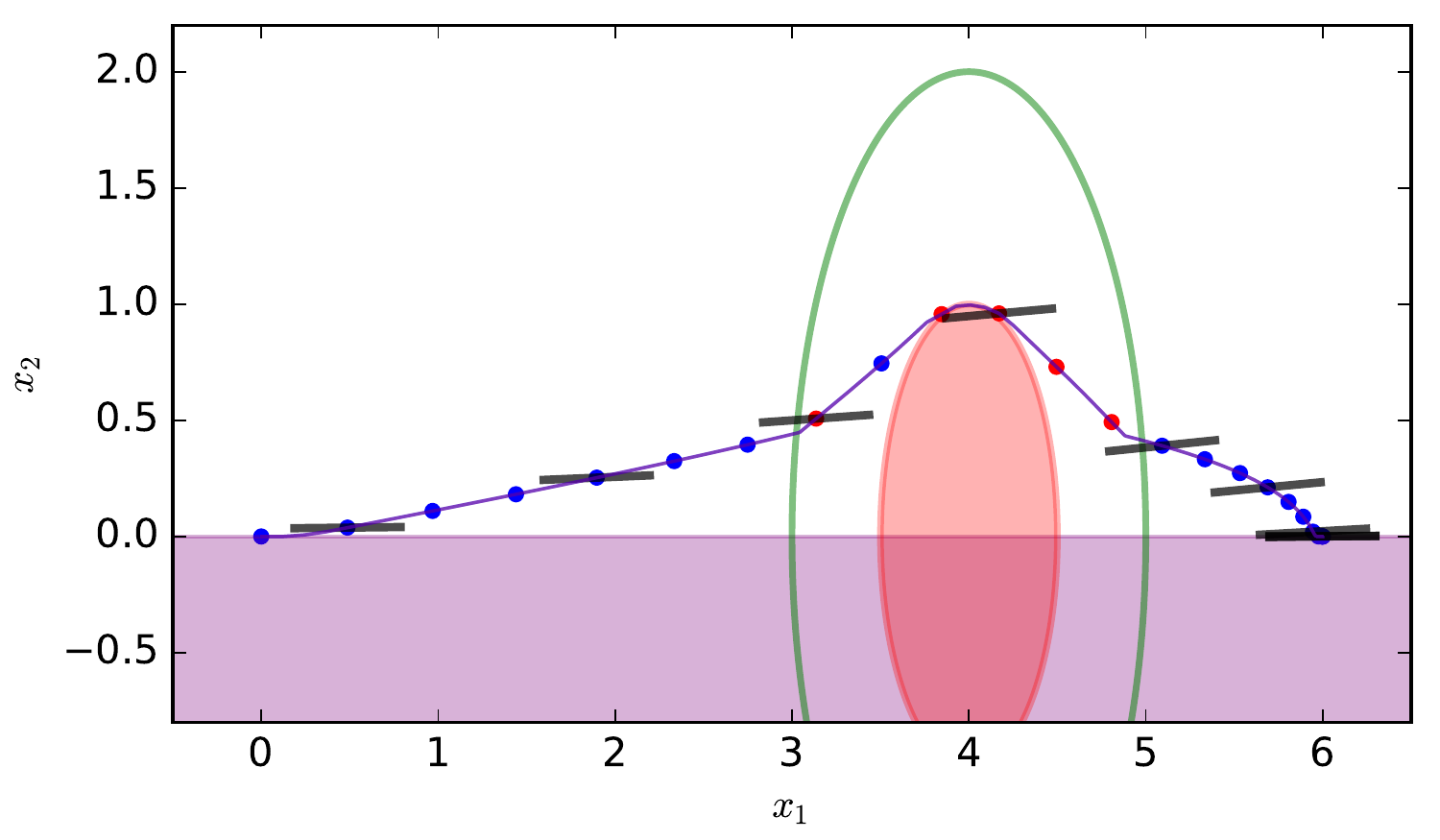}%
  }%
  \hfill%
  \subfloat[%
  Hybrid trajectories using the projection operator ${\cal P}_N\pb{{\cal F}_N(\omega)}$, with $N=6$ (lighter), $N=9$, and $N=12$ (darker).
  Note how the approximation quickly converges to the optimal relaxed trajectory.
  ]{%
    \label{fig:ex3b}%
    \includegraphics[width=.32\linewidth]{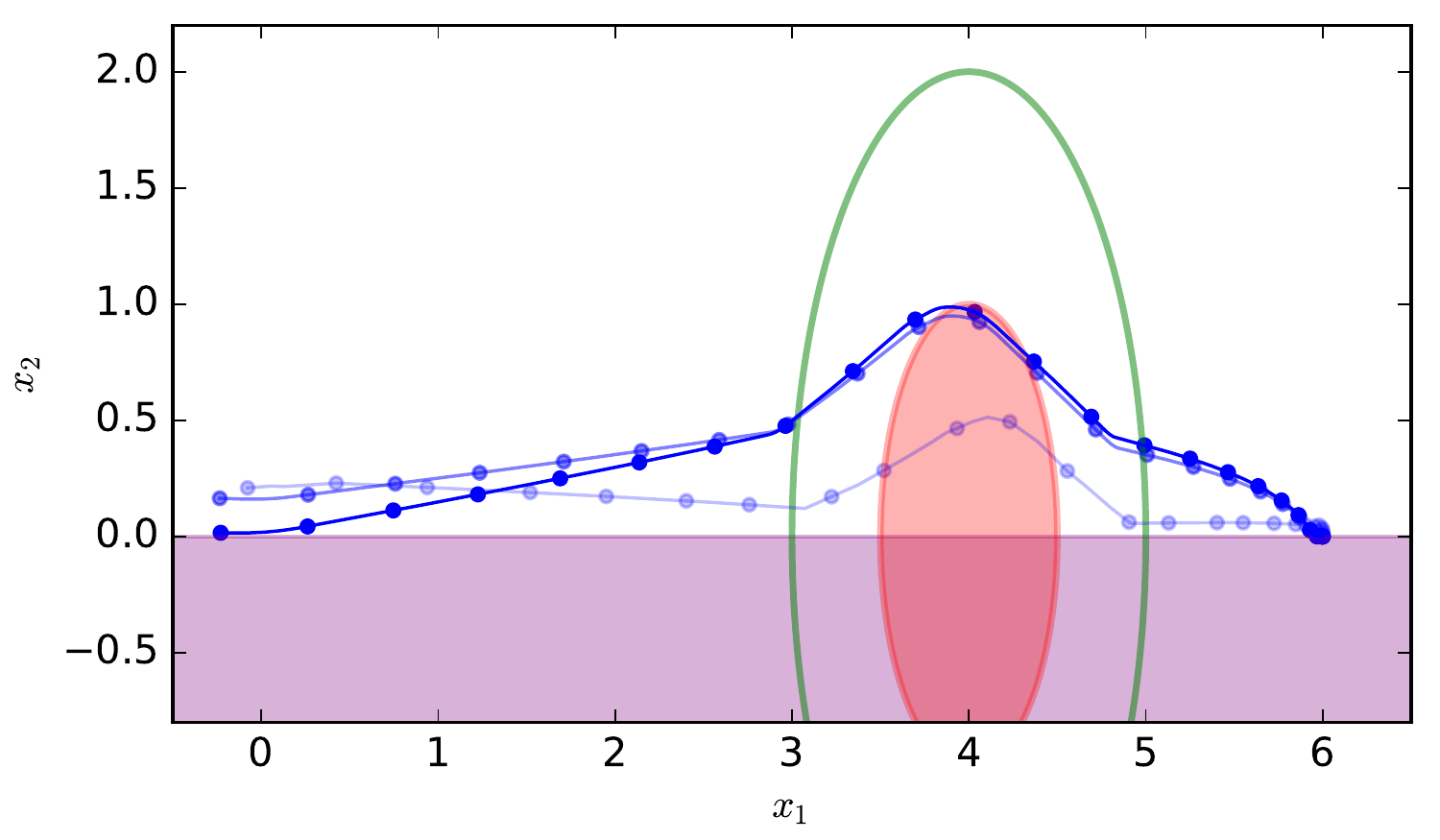}%
  }%
  \hfill%
  \subfloat[%
  Optimal relaxed indexing function $\omega$.%
  ]{%
    \label{fig:ex3c}%
    \includegraphics[width=.32\linewidth]{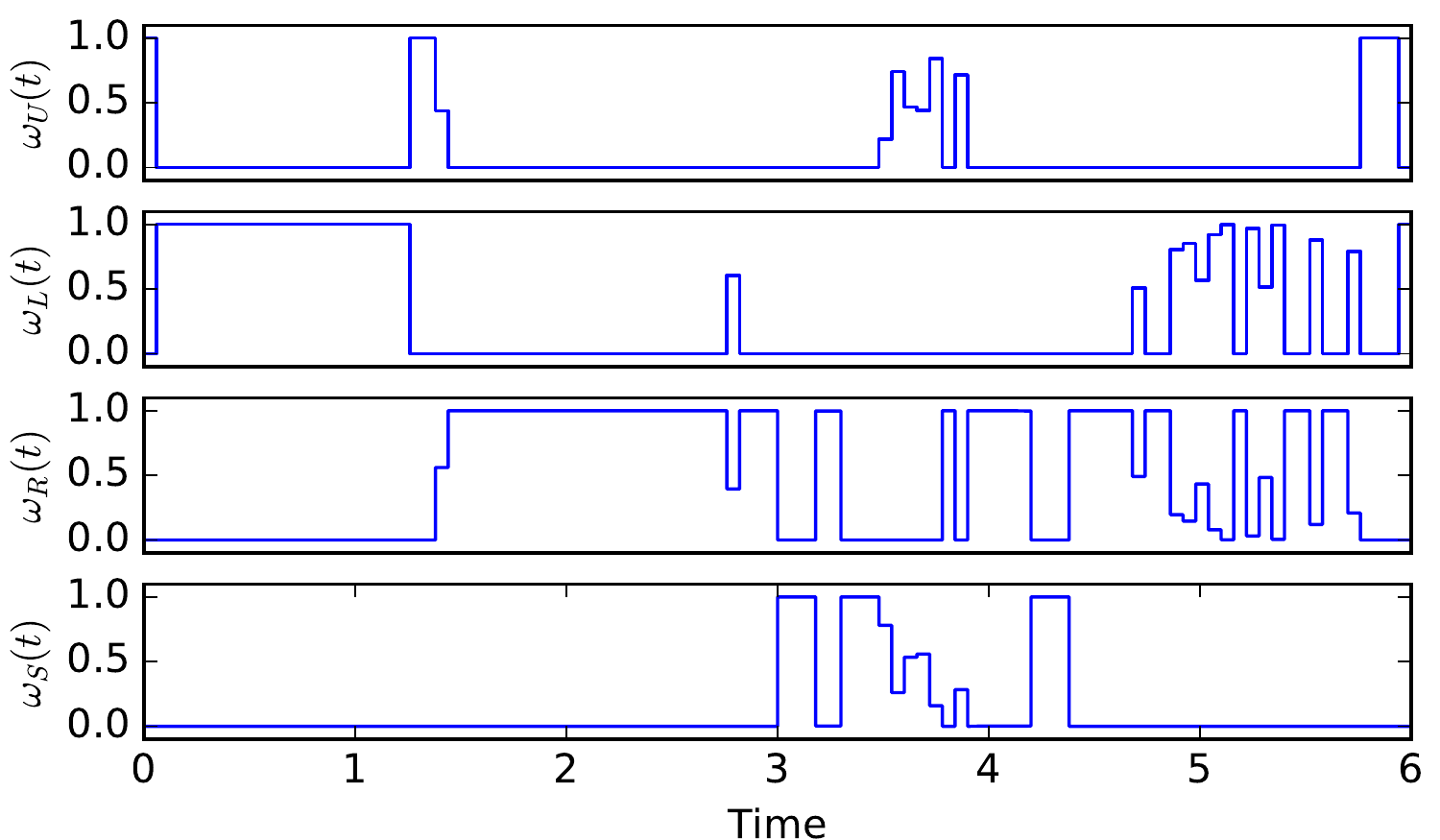}%
  }
  \caption{Results of Example~\ref{sec:ex_qr}.}
\end{figure*}

As shown in Figures~\ref{fig:ex3a} and~\ref{fig:ex3c}, the safety mode fires appropriately to ensure that the helicopter avoids the obstacle, successfully reaching the origin at the end of the simulation.
This example shows that using a hybrid formulation to incorporate auxiliary modes of operation is a suitable strategy to help a control engineer getting the desired behavior in view of unsafe conditions.
In Figure~\ref{fig:ex3b} we show the results of using the operators ${\cal F}_N$ and ${\cal P}_N$ to obtain hybrid indexing functions from the relaxed optimal solution, where $N=12$ produces a very good approximation of the optimal result.

\section{Conclusion}
\label{sec:conclusion}

We have presented a novel mathematical formulation that allows us to numerically solve hybrid optimal control problems for the class of systems whose trajectories are continuous.
Our formulation uses the theory of relaxed controls together with a feedback rule which disables certain discrete modes as a function of the continuous state at any given time.
We have shown that the formulation yields an accurate approximation to the optimal trajectories of a hybrid system in an efficient, and easy to compute, transformation that results in a nonlinear programming problem.

Our future efforts will be focused on extending this formulation to hybrid systems with discontinuous jumps, and to provide a simple computational interface to implement these optimal control problems built as a module of the \emph{OptWrapper} library~\cite{optwrapper}.


\bibliographystyle{IEEEtran}
\bibliography{refs}


\end{document}